\documentclass[12pt]{amsart}

\newtheorem{theorem}{Theorem}[section]
\newtheorem{lemma}[theorem]{Lemma}

\newtheorem{proposition}[theorem]{Proposition}
\newtheorem{remark}[theorem]{Remark}

\newcommand{\ncom}{\newcommand}
\ncom{\ep}{\epsilon}
\ncom{\rar}{\rightarrow}
\ncom{\lrar}{\longrightarrow}
\ncom{\ov}{\overline}
\ncom{\what}{\widehat}

\newcommand{\ignore}[1]{}
\ncom{\m}{\mbox}
\ncom{\sta}{\stackrel}
\ncom{\comx}{{\mathbb C}}
\ncom{\A}{{\mathbb A}}
\ncom{\Z}{{\mathbb Z}}
\ncom{\Q}{{\mathbb Q}}
\ncom{\R}{{\mathbb R}}
\ncom{\G}{{\mathbb G}}
\ncom{\hH}{{\mathbb H}}
\ncom{\al}{\alpha}
\ncom{\p}{{\mathbb P}}
\ncom{\E}{{\mathbb E}}
\ncom{\N}{{\mathbb N}}
\ncom{\K}{{\mathbb K}}
\ncom{\X}{{\mathbb X}}
\ncom{\f}{\frac}
\ncom{\cA}{{\mathcal A}}
\ncom{\cB}{{\mathcal B}}
\ncom{\cD}{{\mathcal D}}
\ncom{\cX}{{\mathcal X}}
\ncom{\cO}{{\mathcal O}}
\ncom{\cW}{{\mathcal W}}
\ncom{\cL}{{\mathcal L}}
\ncom{\cP}{{\mathcal P}}
\ncom{\cH}{{\mathcal H}}
\ncom{\cS}{{\mathcal S}}
\ncom{\cM}{{\mathcal M}}
\ncom{\cC}{{\mathcal C}}
\ncom{\cT}{{\mathcal T}}
\ncom{\cF}{{\mathcal F}}
\ncom{\cN}{{\mathcal N}}
\ncom{\cJ}{{\mathcal J}}
\ncom{\cV}{{\mathcal V}}
\ncom{\cZ}{{\mathcal Z}}
\ncom{\cU}{{\mathcal U}}
\ncom{\cSU}{{\mathcal S \mathcal U}}
\ncom{\cG}{{\mathcal G}}
\ncom{\cQ}{{\mathcal Q}}
\ncom{\cR}{{\mathcal R}}
\ncom{\cY}{{\mathcal Y}}
\ncom{\cE}{{\mathcal E}}
\ncom{\cI}{{\mathcal I}}
\ncom{\mylabel}[1]{{\rm (#1)}\label{#1}}
\ncom{\Hom}{{\textit{Hom}}}
\ncom{\eop}{{\hfill $\Box$}}
\begin{document}
\baselineskip=16pt

%%%%%%%%%%%%%%%%%%%%%%%%%%%%%%%%%%%%%%%%%%%%%%%%%%%%%%%%%%%%%%%%%%%%%%%%%%%%%%%%%%%%%%%%%%%%%%%%%%%%%

\title[Chern-Simons classes]{Chern-Simons classes of flat connections on supermanifolds }

%%%%%%%%%%%%%%%%%%%%%%%%%%%%%%%%%%%%%%%%%%%%%%%%%%%%%%%%%%%%%%%%%%%%%%%%%%%%%%%%%%%%%%%%%%%%%%%%%%%%

\author[J. N. Iyer]{Jaya NN  Iyer}
\address{School of Mathematics, Institute for Advanced Study, 1 Einstein Drive, Princeton NJ 08540 USA.}
\email{jniyer@ias.edu}

\address{The Institute of Mathematical Sciences, CIT
Campus, Taramani, Chennai 600113, India}
\email{jniyer@imsc.res.in}

\author[U.N.Iyer]{Uma N Iyer}
\address{308A, Department of Mathematics and Computer Science,
CP315, Bronx Community College,
University Avenue and West 181 Street,
Bronx, NY 10453.}
\email{uma.iyer@bcc.cuny.edu}

\footnotetext{Mathematics Classification Number: 53C05, 53C07, 53C29 }

\begin{abstract}
In this note we define Chern-Simons classes of a superconnection $D+L$ on a complex supervector bundle $E$ such that $D$ is flat and preserves the grading, and $L$ is an odd endomorphism of $E$ on a supermanifold. As an application we obtain a definition of Chern-Simons classes of a (not necessarily flat) morphism between flat vector bundles on a smooth manifold. 
We extend Reznikov's theorem on triviality of these classes 
when the manifold is a compact K\"ahler manifold or a smooth complex quasi--projective variety, in degrees $>\,1$.
\end{abstract}
\maketitle

%%%%%%%%%%%%%%%%%%%%%%%%%%%%%%%%%%%%%%%%%%%%%%%%%%%%%%%%%%%%%%%%%%%%%%%%%%%%%%%%%%%%%%%%%%

\section{Introduction}

%%%%%%%%%%%%%%%%%%%%%%%%%%%%%%%%%%%%%%%%%%%%%%%%%%%%%%%%%%%%%%%%%%%%%%%%%%%%%%%%%%%%%%%%%%%

Suppose $(X,\cC^{\infty}_X)$ is a $\cC^{\infty}$-differentiable manifold endowed with the structure sheaf
$\cC^{\infty}_X$ of smooth functions.
Let $E$ be a complex $\cC^{\infty}$ vector bundle on $X$ of rank $r$ and equipped with a connection $\nabla$.
The Chern-Weil theory defines the Chern classes
$$
c_i(E,\nabla) \in H^{2i}_{dR}(X,\comx), \m{  for }i=0,1,...,r
$$
in the de Rham cohomology of $X$.
These classes are expressed in terms of the $GL_r$-invariant polynomials evaluated on the curvature form
$\nabla^2$.

Suppose $E$ has a flat connection, i.e., $\nabla^2=0$. Then the de Rham Chern classes are zero. It is significant to define Chern-Simons classes for a flat connection. These are classes in the $\comx/\Z$-cohomology and were defined by Chern-Cheeger-Simons in \cite{Ch-Si}, \cite{Chn-Si}.   

Quillen has pointed out in \cite{Quillen},\cite{Quillen2}, a homomorphism 
$u:E_0\rar E_1$ between vector bundles on a smooth manifold $M$ and inducing an isomorphism over a subset 
$A\subset M$ corresponds to an element in the relative $K$-group $K(M,A)$. A Chern character in the de Rham cohomology of $M$ associated to the homomorphism $u$ is computed in \cite{Quillen} whose class is shown to be equal to the difference $\m{ch}(E_0)- \m{ch}(E_1)$ of the Chern characters. This describes the Chern character of the homomorphism $u$. In fact, we think that it would be good to look at a quiver, i.e., a sequence of homomorphisms between vector bundles 
$$
E_0\rar E_1\rar...\rar E_r
$$
over a smooth manifold and define the Chern character of the sequence in the de Rham cohomology. This will involve a study of $\Z_{r+1}$-graded objects, which we will look in the future. 
Quillens proof involves
regarding $E=E_0\oplus E_1$ as a 
supervector bundle on $M$ and $D$ be any connection preserving the grading and associating an odd endomorphism of $E$, with respect to $u$ and a choice of a metric. 

In this paper, we want to look at a morphism $u$ between flat vector bundles and extend Quillen's
construction and define Chern-Simons classes for the morphism $u$. 
Hence it is relevant to define Chern-Simons classes for flat connections in the setting of supermanifolds, in a more general set-up.

For the definition of \textit{supermanifolds}, due to F.A.Berezin and D. Leites, see \cite{De-Mo} (as well as \cite{BL}, \cite{Le}).
The Chern classes of supervector bundles are defined in \cite{Be-Re}, 
on a supermanifold in the integral cohomology. We note that the usual Chern-Weil theory on smooth manifolds expresses de Rham Chern classes in terms of $GL_n$-invariant polynomials on the curvature form of a connection on a smooth vector bundle. In the supersetting, a study of the invariant polynomials has been carried out
by Sergeev \cite{Se}, following works by Berezin \cite{Be},\cite{Be2} and Kac \cite{Ka}. We do not know if the differential forms defined using the invariant polynomials of Sergeev give the de Rham Chern class of a super vectorbundle, as obtained by Quillen.
In this paper we consider the Chern character defined by Quillen to define the Chern-Simons classes.

Let $(M,\cO_M)$ denote a complex supermanifold and $(M_B,\cC^\infty_M)$ 
denote the underlying $\cC^\infty$-manifold.

With notations as in \cite{De-Mo} or \S \ref{prelim}, we show
\begin{theorem}
Suppose $(M,\cO_M)$ is a complex supermanifold. Let $\cE^{r|s}$ be a complex 
supervector bundle on $(M,\cO_M)$ equipped with a superconnection $\nabla=D+L$ such that $D$ preserves the grading and is flat, and $L$ is an odd endomorphism of $\cE^{r|s}$. 
Then there exists uniquely determined Chern-Simons classes
$$
\widehat{c_i}( \cE^{r|s},\nabla)\, \in \, H^{2i-1}(M_B,\comx/\Z)
$$
for $i>0$. Furthermore, if $M_B$ is a compact K\"ahler manifold or a smooth complex quasi--projective variety, then these classes are 
torsion, in degrees $>\,1$.
\end{theorem}

This can be thought of as an extension of Reznikov's 
fundamental theorem \cite{Reznikov2} on rationality of Chern-Simons classes on compact K\"ahler manifold, in the setting of supermanifolds.   
We also define Chern-Simons classes of a (not necessarily flat) homomorphism $u:E_0\rar E_1$ between flat complex vector bundles, extending Quillen's construction of the de Rham Chern character.
Then we prove a relative Reznikov theorem (see Theorem \ref{RelReznikov}) for the classes of
the morphism $u$.

%-------------------------------------------------------------------------------------

\section{Preliminaries}\label{prelim}
%-------------------------------------------------------------------------------------

We briefly recall the definitions and terminologies from \cite{Le} 
and from the notes by Deligne
and Morgan \cite{De-Mo}.

Let $\cC^\infty$ be the sheaf of $C^\infty$-functions on $\R^p$. The space $\R^{p|q}$ is the topological space $\R^p$, endowed with the sheaf $\cC^\infty[\theta_1,...,\theta_q]$ of supercommutative super $\R$-algebras, freely generated over $\cC^\infty$ by the anticommuting $\theta_1,...,\theta_q$.
The coordinates $t^i$ of $\R^p$ and the $\theta_j$ and all generators of $\cC^\infty$ obtained from them by any automorphism are said to be the coordinates of $\R^{p|q}$. A supermanifold $M$ of dimension $p|q$ is a topological space $M_B$ (or also called as the body manifold with the structure sheaf $\cC^\infty_M$) endowed with a sheaf of super $\R$-algebras which is locally isomorphic to $\R^{p|q}$. The structure sheaf of $M$ is denoted by $\cO_M$.
We denote $p|q$, the real dimension of the supermanifold $M$.

On $M=\R^{p|q}$, the even derivations $\partial/\partial t^i$ and the odd derivations $\partial/\partial\theta^j$ are defined.

\begin{proposition}\cite[2.2.3]{Le}
The $\cO_M$-module of $\R$-linear derivations of $\cO_M$ is free of dimension $p|q$, with basis: the
 $\partial/\partial t^i$ and the $\partial/\partial\theta^j$.
\end{proposition}

Complex supermanifolds are topological spaces endowed with a sheaf of super $\comx$-algebras, locally isomorphic to some $(\comx^p,\cO[\theta^1,...,\theta^q])$. Here $\cO$ is the sheaf of holomorphic
functions on $\comx^p$. As before we denote $p|q$, the complex dimension of the complex supermanifold $M$.

Suppose $R$ be a commutative superalgebra and the standard free module $A^{r|s}$ is the module freely generated by even elements $e_1,...,e_r$ and odd elements $f_1,...,f_s$. 
An automorphism of $A^{r|s}$ is represented by an invertible matrix 
\begin{eqnarray}\label{matrix}
X = \left(  
\begin{array}{cc}
X_1 & X_2 \\
X_3 & X_4
\end{array}
\right) 
\end{eqnarray}
such that the $(r\times r)$--matrix $X_1$ and the $(s\times s)$--matrix $X_4$ have even entries and the $(s\times r)$--matrix $X_3$ and the $(r\times s)$--matrix $X_2$ have odd entries.
The group of all automorphisms of $A^{r|s}$ is denoted by $GL(r,s)$.

Suppose $M$ is a supermanifold and locally it looks like $\R^{p|q}$ as above. 
A \textit{complex supervector bundle} $V$ on $M$ is a fiber bundle $V$ over $M$ with typical fiber $\comx^{r|s}$ and structural group $GL(r,s)$. Alternately, it can be considered as a sheaf of $\cO_M$-supermodules  $\cV$, locally free of rank $r|s$. 

The tangent bundle $\cT_M$ is the $\cO_M$-module of derivations of $\cO_M$ and is a supervector bundle of rank $p|q$. The cotangent bundle $\Omega^1_M$ is the dual of $\cT_M$.
There is a differential $d:\cO_M\lrar \Omega^1_M$, giving rise to the super de Rham complex $\Omega^\bullet_M$ on $M$.

\begin{lemma}{\rm{(Poincar\'e lemma)}}\cite[p.73]{De-Mo}
The complex $\Omega^\bullet_M$ is a resolution of the constant sheaf on the body manifold $M_B$.
\end{lemma} 

In particular, the cohomology of $M_B$ can be computed by the super de Rham complex:
$$
H^*(M_B,\R)=H^*(\Gamma(M,\Omega^\bullet_M)).
$$

We briefly review the group of differential characters and Chern-Simons classes on a smooth manifold $X$.

\subsection{Analytic differential characters on $X$ \cite{Ch-Si}}\label{diff.char}

Let $S_k(X)$ denote the group of $k$-dimensional smooth singular chains on $X$, with integer coefficients. Let $Z_k(X)$ denote the subgroup of cycles.  Let us denote 
$$
S^\bullet(X,\Z):=\m{Hom}_\Z(S_\bullet(X),\Z)
$$
the complex of $\Z$ -valued smooth singular cochains, whose boundary operator is denoted by $\delta$.
The group of smooth differential $k$-forms on $X$ with complex coefficients is denoted by $A^k(X)$
and the subgroup of closed forms by $A^k_{cl}(X)$.
Then $A^\bullet(X)$ is canonically embedded in $S^\bullet(X)$, by integrating forms against the smooth
singular chains. In fact, we have an embedding
$$
i_\Z:A^\bullet(X)\hookrightarrow S^\bullet(X,\comx/\Z).
$$

The group of differential characters of degree $k$ is defined as
$$
\widehat{H^{k}}_\comx(X):=\{(f,\al)\in \m{Hom}_\Z(Z_{k-1}(X),\comx/\Z) \oplus A^k(X): \delta(f)=i_\Z(\al) \m{ and } d\al=0 \}.
$$
There is a canonical and functorial exact sequence:
\begin{equation}\label{diffeqn}
0\lrar H^{k-1}(X,\comx/\Z)\lrar \widehat{H^{k}}_\comx(X)\lrar A^k_{\Z}(X)\lrar 0.
\end{equation}
 
Here $A^k_{\Z}(X):=\m{ker}(A^k_{cl}(X)\lrar H^k(X,\comx/\Z))$.

Similarly, one can define the group of differential characters $\what{H^{k}}_\R(X)$ 
which have $\R/\Z$-coefficients.

%%%%%%%%%%%%%%%%%%%%%%%%%%%%%%%%%%%%%%%%%%%%%%%%%%%%%%%%%%%%%%%%%%%%%%%%%%%%%%%%%%%%%%%%%%

\subsection{Cheeger-Chern-Simons classes}\label{Secondary classes}

Suppose $(E,\theta)$ is a vector bundle with a connection on $X$.
Then the characteristic forms
$$
c_k(E,\theta) \in A^{2k}_{cl}(X,\Z)
$$
for $0\leq k\leq r$=rank $(E)$, are defined using the universal Weil homomorphism \cite{Chn-Si}.
 
The characteristic classes 
$$
\widehat{c_k}(E,\theta) \in \widehat{H^{2k}}_\comx(X)
$$
are defined in \cite{Ch-Si} using a factorization of the universal Weil homomorphism. These classes are functorial lifting of the forms $c_k(E,\theta)$.

Similarly, there are classes
$$
\widehat{c_k}(E,\theta) \in \widehat{H^{2k}}_\R(X).
$$

\begin{remark}\label{lift}
If the forms $c_k(E,\theta)$ are zero, then the classes $\widehat{c_k}(E,\theta)$ lie in the cohomology $ H^{2k-1}(X,\comx/\Z)$. If $(E,\theta)$ is a flat bundle, then $c_k(E,\theta)=0$ and the classes $\widehat{c_k}(E,\theta)$ are called as the Chern-Simons classes of $(E,\theta)$. Notice that the class depends on the choice of $\theta$.
\end{remark}

Beilinson's theory of universal Chern-Simons classes yield classes for a flat connection $(E,\theta)$,
$$
\widehat{c_k}(E,\theta) \in H^{2k-1}(X,\comx/\Z)
$$
which are functorial and additive over exact sequences (see \cite{DHZ} and \cite{Esnault} for another construction).

%----------------------------------------------------------------------------------------

\section{Chern-Simons classes of flat connections on supermanifolds}

%%%%%%%%%%%%%%%%%%%%%%%%%%%%%%%%%%%%%%%%%%%%%%%%%%%%%%%%%%%%%%%%%%%%%%%%%%%%%%%%%%%%%%%%%%%%

Let $(M,\cO_M)$ be a complex supermanifold of dimension $p|q$.
Consider the sheaf of differentials $\Omega^1_M$ 
on $M$ and let $\cE^{r|s}$ be a complex supervector bundle on $M$ of rank $r|s$.

\begin{lemma}\label{vectorbundle}
Given a complex supervector bundle $\cE^{r|s}$ of rank $r|s$ on $M$, there exists 
a direct sum decomposition
$$
E=E_0\oplus E_1
$$
for some complex smooth vector bundles $E_0$ and $E_1$ of rank $r$ and rank $s$ respectively, on the underlying $\cC^\infty$-manifold $M_B$.
\end{lemma} 
\begin{proof}
A rank $r|s$ complex supervector bundle $\cE^{r|s}$ determines two complex smooth bundles $E_0$ and $E_1$ on the underlying smooth manifold $M_B$ of $M$ as follows. One considers the body map
$$
\cO_M\lrar \cC^\infty_M\otimes \comx
$$
which is obtained by forgetting the local anticommuting variables $\theta_j$. Let $\ov{\cE^{r|s}}$ denote the sheaf of (super)sections of $\cE^{r|s}$. Then $\ov{E_{r+s}}:=\ov{\cE^{r|s}}\otimes_{\cO_M} (\cC^\infty_M\otimes \comx)$ is the sheaf of sections of a rank $r+s$ smooth 
complex vector bundle $E_{r+s}$ on the body manifold $M_B$.
Locally, the sheaf $\cE^{r|s}$ is generated by $r$ even elements and $s$ odd elements as a 
$\cO_M=\cC^\infty_M[\theta_1,...,\theta_q]$--module. Hence tensoring with $\cC^\infty_M$ locally
gives a rank $r+s$ free $\cC^\infty_M$--module given by the same generators.
This implies that the complex vector bundle $E_{r+r}$ is of rank $r+s$.
Now, we notice that the structural group of $\cE^{r|s}$ is $GL(r,s)$ and the structural group of the vector bundle $E_{r+s}$ factors via the projection
$$
GL(r,s)\rar GL(r+s).
$$
Using the description of the elements in $GL(r,s)$ in \eqref{matrix}, it follows that
the image under this projection consists of block diagonal matrices of sizes $r\times r$ and $s\times s$. 

This implies that the matrix of the transition functions of $E_{r+s}$ is a block diagonal 
matrix of rank $r$ and rank $s$ which correspond to smooth complex vector bundles $E_0$ and $E_1$ such that $r=\m{rank }E_0$ and $s=\m{rank }E_1$ on $M_B$.
\end{proof}

In view of the above lemma, we may regard a supervector bundle $\cE^{r|s}$ on $M$ as a supervector bundle
$E=E_0\oplus E_1$, on the underlying $\cC^\infty$-manifold $M_B$ where $E_0$ and $E_1$ are $\cC^\infty$-vector bundles on $M_B$.
 
\subsection{Superconnections}
Let $\cE^{r|s}= \,E=\,E_0\oplus E_1$ be a complex supervector bundle on a manifold $M_B$. Let $\Omega(M_B)=\oplus \Omega^p(M_B)$ be the algebra of smooth differential forms with complex coefficients.
Let
$$
\Omega(M_B,E):=\Omega(M_B) \otimes_{\Omega^0(M_B)}\Omega^0(M_B,E).
$$
where $\Omega^0(M_B,E)$ is the space of (super)sections of $E$.

Then $\Omega(M_B,E)$ has a grading with respect to $\Z\times \Z_2$

A \textit{superconnection} $D$ on $\cE^{r|s}$ is an operator on $\Omega(M_B,E)$ of odd degree
satisfying the derivation property
$$
D(\omega\alpha)=(d\omega)\alpha +(-1)^{deg\,\omega}\nabla\alpha.
$$
For example, a connection on $E$ preserving the grading when extended to an operator on $\Omega(M_B,E)$ in the usual way determines a superconnection.

The \textit{curvature} of a superconnection $D$ is the even degree operator $D^2:= D\circ D$ on $\Omega(M_B,E)$.

A superconnection is said to be \textit{flat} if $D^2=0$. We call the pair $(\cE^{r|s},D)$ as a flat complex supervector bundle.

We want to define Chern-Simons classes of $(E,D)$ when $D$ is a flat superconnection.
It is not immediately clear that $D$ induces a flat connection on the individual bundle $E_0$ and $E_1$.

For this purpose, we look at the situation, considered by Quillen \cite{Quillen}.

%%%%%%%%%%%%%%%%%%%%%%%%%%%%%%%%%%%%%%%%%%%%%%%%%%%%%%%%%%%%%%%%%%%%%%%%%%%%%%%%%%%%%%%
\subsection{Quillen's construction}

Suppose $M$ is a supermanifold and $E$ is a complex supervector bundle on $M$. Regard $E=E_0\oplus E_1$ as a complex supervector 
bundle on $M$ in view of Lemma \ref{vectorbundle}, where $E_0$ and $E_1$ are smooth vector bundles on the body manifold $M_B$. Under this identification we omit the suffix $B$ from $M_B$ and without confusion we write $M=M_B$ in the following discussion.

Suppose $E$ is equipped with a superconnection $D$.
From the curvature $D^2$, one can construct differential forms
$$
\m{str}(D^2)^n= \m{str}D^{2n}
$$
in $\Omega(M)^{even}$. Here $\m{str}$ denotes the supertrace.
These are even forms since the supertrace preserves the grading. 

We have,
\begin{theorem}\label{quillen}
The form $\rm{str}D^{2n}$ is closed, and its de Rham cohomology class is independent of the choice of superconnection $D$.
\end{theorem}
\begin{proof}
See \cite[Theorem, p.91]{Quillen}.
\end{proof}

Quillen described the (super) Chern character of $E$ in terms of the usual Chern characters of $E_0$ and $E_1$ in the following situation.

Regard $E=E_0\oplus E_1$ as a complex supervector bundle and $D_0=D^0+D^1$ be a connection on $E$ preserving the grading. Let $L$ be an odd degree endomorphism of $E$ and write $D_t:=D+t.L$ where $t$
is a parameter.

\begin{proposition}$\rm{(Quillen)}$\cite{Quillen}\label{Quillen}
Replacing $L$ by $t.L$, where $t$ is a parameter, one obtains a family of forms 
\begin{equation}\label{traceform}
\rm{str}\,e^{(D+tL)^2}\,=\,\rm{str}\,e^{r^2L^2+ t[D,L]+ D^2}
\end{equation}
all of which represent the Chern character $\rm{ch}(E_0)-\rm{ch}(E_1)$ in the de Rham cohomology of $M$.

\end{proposition}
\eop

The referee has pointed out that the above computations on a supermanifold produces pseudodifferential forms. For our purpose, it suffices to note that the trace form in \eqref{traceform} defines a closed differential form whose de Rham class is independent of the superconnection \cite[Theorem p.91]{Quillen}.

 In this situation we want to define uniquely determined
 Chern-Simons classes of $(E,D_t)$ which is independent of $t$ and if $D^0, D^1$ are flat connections.
For this purpose, we look at the \textit{Character diagram} of Simons and Sullivan (see \cite{Sm-Su}):

\begin{equation}\label{Characterdiagram}
\begin{array}{ccccccc}
0 & & & & & & 0 \\
 &&&&&& \\
  & \searrow & & & & \nearrow &  \\
 &&&&&& \\
  &          & H^{k-1}(\R/\Z) & \sta{-B}{\lrar} & H^k(\Z) & & \\
&&&&&& \\
\searrow & \,\,\,\,\,\, \alpha{\nearrow} &\,\,\,\,\,\,\,\, \searrow i_1 & \,\,\,\,\,\,\,\,\,\,\,\delta_2\nearrow &\,\,\,\,\,\,\,\,\,\, \searrow r  & \,\,\,\,\,\,\,\,\,\nearrow \\
&&&&&& \\
& H^{k-1}(\R) & &\what{H^k}_\R(M)  && H^k(\R) & \\
&&&&&& \\
\nearrow &\,\,\,\,\,\,\,\,\, \searrow \beta &\,\,\,\,\,\,\,\,i_2 \nearrow &\,\,\,\,\,\,\,\,\,\,\,\,\, \searrow \delta_1 &\,\,\,\,\,\,\,\,\,\,\,\,\,\,s \nearrow &\,\,\,\,\,\,\,\,\,\,\,\,\, \searrow \\
&&&&&&\\
         &          & A^{k-1}/A^{k-1}_{\Z} & \sta{d}{\lrar}& A^k_{\Z} && \\
&&&&&& \\
 & \nearrow &&&& \searrow &   \\
&&&&&&\\
0 &&&&&& 0
\end{array}
\end{equation}

The diagonal sequences are exact and $(\alpha,B,r)$ is the Bockstein long exact sequence associated to the coefficient sequence $\Z\rar \R\rar \R/\Z$. Also $(\beta,d,s)$ is another long exact sequence in which $\beta$ and $s$ are defined via the de Rham theorem.
(A similar diagram holds by replacing $\R$ with $\comx$).

\begin{lemma}\label{unique}
Suppose $(F,\nabla)$ is a smooth flat connection on a manifold $M$. Then the Chern class
$c_i(F)\in H^{2i}(M,\Z)$ vanishes in $H^{2i}(M,\R)$. There is a unique lifting 
$\what{c_i}(F,\nabla)\in H^{2i-1}(M,\R/\Z)$ of the  integral class $c_i(F)$, for $i>0$.
\end{lemma}
\begin{proof}
The first assertion is by the Chern-Weil theory. The second assertion is by the Chern-Simons-Cheeger construction \cite{Ch-Si}.
\end{proof}

Consider the total Chern class
$$
c(F)= 1+c_1(F)+...+c_f(F)
$$
and the total Segre class
$$
s(F)=1+s_1(F)+...+s_f(F).
$$
Then we have the relations
\begin{equation}\label{segre}
s(F)=\f{1}{c(F)}
\end{equation}
and 
\begin{equation}\label{sum}
c(F-G)=c(F).s(G)
\end{equation}
where $G$ is any other vector bundle.
These relations also hold if we replace the classes $c_i(F)$ by $\what{c_i}(F,\nabla)$ and $s_i(F)$ by $\what{s_i}(F,\nabla)$ which are defined by the relation \eqref{segre}.
See \cite[p.64-65]{Ch-Si}.

To define a canonical lifting in $\what{H^*_\R}(M)$ of the trace form \eqref{traceform} associated to the superconnection $(E,D_t)$, one would need a notion of universal superconnection analogous to the universal connections defined by Narasimhan and Ramanan \cite{Narasimhan}, \cite{Narasimhan2}, which we may look in the future.

Nonetheless, we consider the family of superconnection $D_t=D+t.L$ as above. Consider the product manifold $\R\times M$ and the pullback $pr_2^*E$ of the bundle $E$. This bundle is equipped with a superconnection
$$
\bar{D}:= dt\partial_t+D'
$$
whose restriction to $\{t\}\times M$ is $D_t$. In terms of local trivialization of $E=M\times V$ we can describe
$\bar{D}$, $D'$ as follows. Write $D_t=d_M+\theta_t$, where $\theta_t$ is a family of one-forms on $M$ with values in $\m{End} V$ and let $\theta$ be the form on $\R\times M$ not involving $dt$ and having the restriction $\theta_t$ on $\{t\}\times M$. Then $D'=d_M+\theta$ and 
$$
\bar{D}=dt\partial_t+D'=(dt\partial_t+d_M)+\theta=d_{\R\times M} +\theta.
$$
See also \cite[p.91]{Quillen}.

By the homotopy property of de Rham cohomology, it follows that the class of $\m{str}D^{2n}_t$ in $H^{2n}(M,\R)$ is independent of $t$. 

\begin{proposition}\label{ChernSimons}
Suppose the superconnection $D_0=D^0\oplus D^1$ on the supervector bundle $E=E_0\oplus E_1$ preserves the grading and the individual smooth connections $D^0$ and $D^1$ are flat connections on $E_0$ and $E_1$ respectively. Then $D_t=D_0+t.L$ is a superconnection on $E$ where $L$ is an odd endomorphism of $E$. 
Then there is a uniquely determined class $\what{c_n}(E,D_t) \in H^{2n-1}(M,\R/\Z)$, independent of $t$. Moreover this class is equal to 
$$
\what{c_n}(E,D_0+L)=\sum_{p+q=n}\what{c_p}(E_0,D^0).\what{s_q}(E_1,D^1)
$$
\end{proposition}  
\begin{proof}
We notice that the trace forms are integral valued. This implies that the Chern class associated
to these forms lies in the integral cohomology $H^{2n}(M,\Z)$ which is independent of $t$ in
$H^{2n}(M,\R)$, by Quillen's Theorem \ref{quillen}. This determines a class in $H^{2n}(M,\Z)$ independent of $t$. But this class vanishes in $H^{2n}(M,\R)$ since $D_0$ has components $D^0$ and $D^1$ which are flat connections, hence $D_0^2=0$. Using the Bockstein operator in \eqref{Characterdiagram}, we conclude that there is a class
$$
\what{c_n}(E,D_t) \,\in\,H^{2n-1}(M,\R/\Z)
$$
which is independent of $t$ and we denote this class by $\what{c_n}(E,D_0+L)$.
We get a uniquely determined class $\what{c_n}(E,D_0+L)=\what{c_n}(E, D_0)\in H^{2n-1}(M,\R/\Z)$, as follows.

To get an expression for this class, we notice that by Quillen's result Proposition \ref{Quillen} the (super) Chern character form of $E$ is the difference $\m{ch}(E_0)-\m{ch}(E_1)$ in the integral cohomology. In particular we want to lift the integral Chern classes of $E_0-E_1$ in the $\R/\Z$-cohomology.
The relations in \eqref{segre} and \eqref{sum} give the formula
$$
\what{c_n}(E,D_0+L):=\sum_{p+q=n}\what{c_p}(E_0,D^0).\what{s_q}(E_1,D^1).
$$
The uniqueness of $\what{c_n}(E,D_0+L)$ follows from the uniqueness of the 
Chern-Simons classes $\what{c_p}(E_0,D^0)$ and $\what{c_q}(E_1,D^1)$, see Lemma \ref{unique}.
This concludes the lemma.
\end{proof}

\begin{remark}\label{sclasses}
All the above constructions follow verbatim by replacing $\R/\Z$-coefficients with $\comx/\Z$-coefficients. We call the resulting classes $\what{c_n}(E,D_0+L)\in H^{2n-1}(M,\comx/\Z)$ as the Chern-Simons classes of the superconnection $D_0+L$ (or $D_0+t.L$ for a parameter $t$).
\end{remark}

%%%%%%%%%%%%%%%%%%%%%%%%%%%%%%%%%%%%%%%%%%%%%%%%%%%%%%%%%%%%%%%%%%%%%%%%%%%%%%%%%%%%%%%
\subsection{Chern Simons classes for a morphism between flat connections}

 Consider a homomorphism $u:E_0\rar E_1$ between complex vector bundles on a smooth manifold $M$. Then $u$ determines a class in the $K$-group $K(M)$.

Let
\begin{eqnarray}\label{L}
L=i \left(  
\begin{array}{cc}
0 & u^* \\
u & 0
\end{array}
\right) .
\end{eqnarray}
Here $u^*$ is the adjoint of $u$ relative to a given metric (see \cite{Quillen}). 
Regard $E=E_0\oplus E_1$ as a complex supervector bundle on $M$ and $D_0=D^0+D^1$ be a superconnection on $E$ preserving the grading. Then $L$ is an odd degree endomorphism of $E$ and as shown in \cite{Quillen}, $D+L$ is a superconnection and its curvature form $F=(D+L)^2$ is an even form with values in $\m{End}E$.

\begin{lemma}\label{relative}
Suppose $(E_0,D^0)$ and $(E_1,D^1)$ are flat connections and $u$ and $L$ are as above. 
Then we can define the Chern-Simons classes of the morphism $u$ (which need not be a flat morphism) 
in the $\comx/\Z$-cohomology of $M$ by setting
\begin{equation}\label{CSmorphism}
\what{c_n}(u)\,:=\,\what{c_n}(E,D_0+L)  \in\,H^{2n-1}(M,\comx/\Z)
\end{equation}
for $n\geq 1$, where $\what{c_n}(E,D_0+L)$ are defined in Lemma \ref{ChernSimons}.
\end{lemma}
\begin{proof}
The assumptions of Proposition \ref{ChernSimons} are fulfilled and we obtain uniquely defined classes
$$
\what{c_n}(u)\,:=\, \what{c_n}(E,D_0+L)  \in\,H^{2n-1}(M,\comx/\Z)
$$
for $n\geq 1$.
\end{proof}

\begin{theorem}$\rm{(Relative \,Reznikov\,theorem)}$\label{RelReznikov}
Suppose $u:E_0\rar E_1$ is a (not necessarily flat) homomorphism between flat complex vector bundles $(E_0,D^0)$ and $(E_1,D^1)$ on a compact K\"ahler manifold $M$ or a smooth complex quasi--projective variety $M$. Then the classes 
$$
\what{c_i}(u) \in H^{2i-1}(M,\comx/\Q)
$$ 
are zero, for $i\geq 2$. 
\end{theorem}
\begin{proof}
By Proposition \ref{ChernSimons}, Remark \ref{sclasses} and Lemma \ref{relative} we have the explicit expression of the class 
$$
\what{c_n}(u)=\sum_{p+q=n} \what{c_p}(E_0,D^0).\what{s_q}(E_1,D^1).
$$
When $M$ is a compact K\"ahler manifold then Reznikov's theorem \cite{Reznikov2} says that
$$
\what{c_n}(E_0,D^0),\, \what{c_n}(E_1,D^1) \in H^{2n-1}(M,\comx/\Z)
$$
are torsion, if $n\geq 2$. A similar result is true if $M$ is a smooth complex quasi--projective variety, by \cite{Iy-Si}. Since the classes $\what{s_q}$ are expressed in terms of $\what{c_i}$ for $i\leq q$, the assertion follows.
This proves the theorem
\end{proof}

{\Small
\textbf{Acknowledgements}:
We thank D. Leites for his comments and pointing out useful 
references to us, especially \cite{Quillen}, \cite{Quillen2}. His comments prompted us to improve
the results and understanding of this subject.
 The first named author is supported by the 
National Science Foundation (NSF) under agreement No. DMS-0111298.  
}

%%%%%%%%%%%%%%%%%%%%%%%%%%%%%%%%%%%%%%%%%%%%%%%%%%%%%%%%%%%%%%%%%%%%%%%%%%%%%%%%%%%%%%%%%%%%%%%%%%%%

\end{document}